\documentclass[reqno,12pt]{amsart} 
\usepackage[top=1.5in,right=1.125in,left=1.125in,bottom=1.5in]{geometry}
\usepackage{amssymb}
\usepackage{color}
\usepackage{tikz}
\usetikzlibrary{positioning,decorations.pathmorphing}
\usetikzlibrary{positioning,decorations.pathmorphing}
\definecolor{red}{rgb}{0.7,0,0}
\definecolor{grey}{RGB}{112,112,112}
\definecolor{blue}{RGB}{034,113,179}
\usepackage[colorlinks=true,citecolor=blue,linkcolor=red,urlcolor=grey,backref=page]{hyperref}
\newcommand{\koniec}{\begin{flushright}  $\Box $ \end{flushright}}
\newtheorem{theo}{Theorem}[section] 
\newtheorem{prop}[theo]{Proposition}  
\newtheorem{lemma}[theo]{Lemma}
\newtheorem{defi}[theo]{Definition}

\theoremstyle{remark}

\newcounter{mnotecount}[section]

\renewcommand{\themnotecount}{\thesection.\arabic{mnotecount}}

\newcommand{\mnote}[1]
{\protect{\stepcounter{mnotecount}}$^{\mbox{\footnotesize
$
\bullet$\themnotecount}}$ \marginpar{
\raggedright\tiny\em
$\!\!\!\!\!\!\,\bullet$\themnotecount: #1} }
\newcommand{\md}[1]{\mnote{{\bf MD:}#1}}

\newcommand{\hook}{{\setlength{\unitlength}{11pt}   
                   \begin{picture}(.833,.8)
                   \put(.15,.08){\line(1,0){.35}}
                   \put(.5,.08){\line(0,1){.5}}
                   \end{picture}}}
\newcommand{\CP}{\mathbb{CP}}
\newcommand{\II}{\mathbb{I}}
\newcommand{\C}{\mathbb{C}}

\newcommand{\F}{\mathbb{F}}
\newcommand{\DD}{\mathbb{D}}
\newcommand{\PP}{\mathbb{P}}
\newcommand{\RP}{\mathbb{RP}}
\newcommand{\R}{\mathbb{R}}
\newcommand{\X}{\mathcal{X}}
\def\p{\partial}
\def\be{\begin{equation}}
\def\D{\mathcal{D}}
\def\OO{\mathcal{O}}
\def\tZ{\widetilde{Z}}
\def\tW{\widetilde{W}}

\def\dim{\mbox{dim}}
\def\rk{\mbox{rank}}
\def\t{\tilde}
\def\id{\mbox{Id}}
\def\ee{\end{equation}}

\def\bea{\begin{eqnarray}}
\def\eea{\end{eqnarray}}
\newcommand{\spp}{\mathbb{S}}
\numberwithin{equation}{section}
\begin{document} \date{September 15, 2021}
\title{Null K\"ahler geometry and isomonodromic deformations}
\author{Maciej Dunajski}
\address{Department of Applied Mathematics and Theoretical Physics\\ 
University of Cambridge\\ Wilberforce Road, Cambridge CB3 0WA, UK.}
\email{m.dunajski@damtp.cam.ac.uk}
\maketitle
\begin{center}
{\em Dedicated to Jenya Ferapontov on the occasion of his 60th birthday.}
\end{center}
\begin{abstract}
We construct the normal forms of null--K\"ahler metrics: pseudo--Riemannian metrics admitting a compatible  parallel nilpotent endomorphism of the 
tangent bundle. Such metrics are examples of non--Riemannian holonomy 
reduction, and (in the complexified setting) appear on the space of Bridgeland stability conditions on a Calabi-Yau threefold.

Using twistor methods we show that, in dimension four - where there is a connection with 
dispersionless integrability - the cohomogeneity--one anti--self--dual null--K\"ahler metrics are generically characterised by solutions
to Painlev\'e I or Painlev\'e II ODEs.
\end{abstract}
\section{Introduction}
A null--K\"ahler structure on a manifold $\X$ of real dimension $4n$ 
consists of a pseudo--Riemannian
metric $g$ of signature $(2n, 2n)$, together with a rank $2n$ endomorphism $N$ of $T\X$ which, for all vector fields   $X, Y$,  satisfies
\[
g(X, NY)+g(NX, Y)=0, \quad N^2=0
\]
and is parallel with respect to the Levi--Civita connection of $g$. 

There are three reasons for considering such structures.
Firstly, they provide an example of a pseudo-Riemannian holonomy reduction with no Riemannan analogue \cite{kath, BI, bryant, GL}
(in the sense that a null--K\"ahler metric can not be analyticaly continued to Riemannian signature). Secondly, there exists a natural fibration of $\X$ over a symplectic
manifold of dimension $2n$, such that the pull--back of the symplectic form to $\X$ 
agrees with the fundamental form $\Omega$ defined by $\Omega(X, Y)=g(NX, Y)$. This structure, albeit in the complexified setup and under the additional condition that $g$ is hyper--K\"ahler, underlies the Bridgeland approach to stability conditions on a
three-dimensional Calabi-Yau triangulated category
\cite{Bridgeland1, Bridgeland2, Bridgeland3}. 
Finally, in dimension four and under the additional anti--self--duality
assumption, null--K\"ahler structures are characterised by solutions to a dispersionless integrable system \cite{D02}.

In the next section we shall introduce null Hermitian structures on
vector spaces, and in \S\ref{section3} we put these structures on manifolds.
Our main result in \S\ref{section3} is the local normal form of the null--K\"ahler condition
\begin{theo}
\label{theorem1}
Let $(\X, g, N)$ be a $4n$--dimensional null--K\"ahler manifold. 
There exist a local coordinate system
$(x^i, y^i), i=1, \dots, 2n$ and a function $\Theta:\X\rightarrow \R$ such that
\begin{eqnarray}
\label{nk_form}
g&=&\frac{1}{2}\sum_{i, j} \omega_{ij}(dy^i\otimes dx^j+dx^j\otimes dy^i) 
+\frac{\p^2 \Theta}{\p y^i \p y^j}(dx^i \otimes dx^j +dx^j \otimes dx^i),\\
 \quad
N&=&\sum_{i} dx^i\otimes \frac{\p}{\p y^i}, \quad\mbox{where}\quad 
\omega_{ij}=\left(\begin{array}{cc}
0&\II_n\\
-\II_n&0
\end{array}
\right).\nonumber
\end{eqnarray}
Conversely, the structure (\ref{nk_form}) is null--K\"ahler for any function $\Theta=\Theta(x^i, y^i)$.
\end{theo}
Thus,  in the real analytic category,  a null--K\"ahler manifold depends on one arbitrary function  of $4n$ variables. In \S\ref{section32} and \S\ref{section33} we shall list
systems of PDEs satisfied by this function if the metric is Einstein or (complexified) hyper--K\"ahler together with examples.

In \S\ref{section4} we focus on oriented four--dimensional null--K\"ahler structures, with the natural choice of orientation which makes the fundamental form $\Omega$ self--dual.
While null--K\"ahler metrics with self--dual Weyl tensor can be found explicitly, the 
anti--self--duality condition on  Weyl tensor corresponds to  
solutions to a 4th order dispersionless integrable system \cite{D02}. 
A problem of 
finding Ricci--flat null reduces to a non--integrable hyper--heavenly equation of Plebanski and Robinson \cite{pleban_robinson}.

In \S\ref{section5} we establish the main result of the paper.
Imposing the invariance
under the isometric action of $SL(2)$ on anti--self--dual null--K\"ahler 
structures leads to an ODE. 
By exploring the connection between the twistor distribution and the isomonodromic Lax pair we shall show that this ODE is either completely solvable, or
transforms to Painlev\'e I or Painlev\'e II.
\begin{theo}
\label{theorem2}
Let $(\X, g, N)$ be an anti--self--dual null--K\"ahler four--manifold with an isometric action 
of $SL(2)$ with three--dimensional orbits, and preserving
the endomorphism $N$. Then either the metric $g$ is conformal to a Ricci--flat metric, or it  can be put in the form
\be
\label{metric_form}
g=\sum_{\alpha, \beta=1}^3\gamma_{\alpha\beta}(t) \sigma^{\alpha}\otimes\sigma^{\beta}
+\sum_{\alpha=1}^3 n_\alpha(t) (\sigma^{\alpha} \otimes dt+dt \otimes \sigma^{\alpha}),
\ee
where the function $t:\X\rightarrow \R$ parametrises the orbits of $SL(2)$, and is constant
on each orbit,  $(\sigma^1, \sigma^2, \sigma^3)$ are left--invariant one--forms on $SL(2)$ which satisfy (\ref{one_formss}),
$\gamma$ is a symmetric 3 by 3 matrix and $n$ is a vector with components 
given by (\ref{solv_form}), or
depending on
solutions of Painlev\'e I or Painlev\'e II as in (\ref{piform}) 
and (\ref{piiform}).
\end{theo}
Anti--self--dual $SL(2)$ invariant Ricci--flat metrics in neutral signature are all known \cite{DT_k}, so the novelty in Theorem \ref{theorem2} lies in the apperance of Painlev\'e equations in the conformal structures with no Ricci--flat representatives.
The connection between the first two Painlev\'e transcendents, and anti--self--dual null--K\"ahler
structures has twistorial origins: the additional structure on the twistor space ${\mathcal Y}$ of $\X$ which corresponds to the endomorophism $N$ is a holomorphic section of $\kappa^{-1/4}$, where $\kappa$ is the holomorphic canonical bundle of ${\mathcal Y}$. The isometric $SL(2)$ action
on ${\mathcal Y}$ also gives rise to such section if the holomorphic vector fields generating this action are linearly dependent at one point, to order four, on each twistor lines. This corresponds
to the isomonodromic problem underlying Painlev\'e I and II.
\subsection*{Acknowledgments} 
I am grateful to
Centro de Investigacion y de Estudios
Avanzados in Mexico for the hospitality 19 years ago, when most of the results in \S\ref{section3} were obtained. My interest in the subject
has been revived after I came across the works of Bridgeland 
\cite{Bridgeland1, Bridgeland2}, and Bridgeland and Strachan \cite{Bridgeland3}, where a null--K\"ahler structure implicitly arises on the space of stability conditions on a 
Calabi--Yau three--fold tirangulated category.
  My work has been partially supported 
by STFC consolidated grants ST/P000681/1, ST/T000694/1.
\section{Algebraic preliminaries}
\label{section2}
Let $V$ be a vector space over field $\F$, where $\F$ is $\R$ or $\C$.
\begin{defi} A null structure on an even--dimensional 
vector space $V$ is an endomorphism $N$ of $V$ such that
\be
\label{eq1}
N^2=0,
\ee
and $\rk(N)=\frac{1}{2}\dim(V)$.
\end{defi}
For reasons to become clear later (see \S\ref{section3}) we shall 
chose $\dim{(V)}=4n$, where $n$ is an integer.
 The kernel of $N$ is a $2n$--dimensional subspace of $V$, and
any basis of this kernel can be extended to a basis of $V$. A convenient way to do it is to pick $2n$ linearly independent vectors $(X_1, X_2, \dots, X_{2n})$
not in $\mbox{Ker}(N)$, and use 
\[
X_1, \dots, X_{2n}, N(X_1), \dots, N(X_{2n})
\]
as a basis of $V$. We shall call this basis adapted to $N$. The matrix of $N$ with respect to an adapted basis is
\be
\label{eq2}
N_0=\left(\begin{array}{cc}
0&\II_{2n}\\
0&0
\end{array}
\right),
\ee
where $\II_{2n}$ is the $2n\times 2n $ identity matrix.

Let $N$ and $N'$ be two null structures on vector spaces $V$ and
$V'$ respectively. A linear map $\phi:V\rightarrow V'$ is called
null linear if
\[
N'\circ \phi=\phi \circ N.
\]
The sub--group ${\mathcal N} (V)\subset GL(V)$ 
of null--linear maps consists of  matrices which commute with the matrix
$N_0$. These matrices are of the form
\[
\left(\begin{array}{cc}
A&B\\
0&A
\end{array}
\right), 
\]
where $A$ and $B$ is an arbitrary $2n\times 2n$ matrix over $\F$, and $A$ is 
invertible.
\begin{prop}
There is a one-to-one correspondence between the set of null
structures on $V$, and the homogeneous space 
$GL(V)/ {\mathcal N}(V)$.
\end{prop}
\noindent
{\bf Proof.}
Consider a $GL(V)$ action of a space of null structures given by
\be
\label{Saction}
N\longrightarrow \phi N\phi^{-1}, \qquad \phi\in GL(V).
\ee
Let $N$ and $N'$ be two null structures of $V$
with adapted bases $(X_i, N(X_i)$ and $({X'}_i, N({X'}_i))$ respectively, where
$i=1,\dots, 2n.$ Define the element $\phi\in GL(V)$ by
\[
\phi(X_i)={X'}_i, \quad \phi(N(X_i))=N'({X'}_i).
\] Therefore $N'=\phi N\phi^{-1}$ and the action
(\ref{Saction}) is transitive. The isotropy subgroup is
of this action is ${\mathcal N}(V)$, as $N_0=\phi N_0\phi^{-1}$ iff $\phi\in
{\mathcal N}(V)$.
\koniec
Our elementary discussion has so far followed the standard
treatment of complex structures \cite{Kob_N}, except the endomorphism $N$ squares to $0$
rather than $-\II_{4n}$. The next step is to introduce the analog of Hermitian inner 
products\footnote{
The natural next step in the theory of complex structures is to introduce a
complexification, where multiplication by complex numbers is given in terms of the complex structure $J$ by
$(a+ib)X=a+b J(X)$, where $X\in V, a, b \in \R$ and $i^2=-1$. Pursuing this analogy for  null structures
leads to dual numbers in place of complex numbers.
The algebra $\DD$ of dual numbers consists of elements of  the form
\[
a+\epsilon\; b, \quad\mbox{where}\quad a, b\in \R\, \quad\mbox{and}\quad
\epsilon^2=0.
\]
The dual numbers can be added, and multiplied according to
\[
(a_1+\epsilon\; b_1)(a_2+\epsilon\; b_2)=a_1 a_2+\epsilon\;(a_1b_2+b_1 a_2), \quad
(a_1+\epsilon \;b_1)+(a_2+\epsilon \;b_2)=a_1+a_2+\epsilon\; (b_1+b_2),
\]
but $\DD$ is not a division algebra, as elements of the form $\epsilon\; b$ do not have inverses. The {\it
infinitesimal} dual number $\epsilon$ underlies non--standard analysis, as it gives a framework to
distinguish between real numbers like $1$, and $0.999\dots$ which are regarded as equal in ordinary
analysis.

A real vector space $V$ with a null structure $N$ can be turned into a vector space $V^{\DD}$ over $\DD$
by defining
\[
(a+\epsilon b) X=aX+b N(X),  \quad\mbox{where}\quad X\in V\quad \mbox{and}\quad a, b \in \R.
\]
In what follows we shall not explore any further connection between null structures and dual numbers,  but
will instead focus on null structures on curved manifolds.}
\begin{defi}
\label{defi2}
A non--degenerate symmetric bi--linear inner product $g:V\times V\rightarrow \F$
is called null--Hermitian if
\be
\label{compat}
g(NX, Y)=-g(X, NY),
\ee
for all $X, Y \in V$, where $N$ is a null structure of $V$. 
\end{defi}
If $\F=\R$, then the signature of a null--Hermitian inner product is $(2n, 2n)$, also called
split, neutral or Kleinian. The definition \ref{defi2} also implies
$
g(X, NX)=0
$ 
and that $\mbox{Ker}(N)$ is a totally isotropic subspace of $V$.

To each null--Hermitian inner product we associate a skew--symmetric bi--linear map
$\Omega\in \Lambda^2(V^*)$ defined by
\be
\label{hermitean_form}
\Omega(X, Y)=g(NX, Y).
\ee
Therefore $\Omega$ vanishes on $\mbox{Ker}(N)$, and it equips the $2n$--dimensional 
quotient vector space $V/\mbox{Ker}(N)$ with a symplectic structure. In a basis
adapted to $N$ the inner product $g$ and the skew--form $\Omega$ are represented by
\[
g_0=\left(\begin{array}{cc}
0&\omega\\
\omega^T&0
\end{array}
\right),\qquad \Omega_0=
\left(\begin{array}{cc}
\omega&0\\
0&0
\end{array}
\right),\qquad\mbox{where}\quad\omega=
\left(\begin{array}{cc}
0&\II_n\\
-\II_n&0
\end{array}
\right).
\]
\subsubsection{Example}
\label{example1}
Let $W$ be a $2n$--dimensional symplectic vector space with a symplectic form $\omega$. The $4n$--dimensional space
$V=W\oplus W$ carries a null--Hermitian structure defined by
\[
N(x, y)=(0, x), \quad g((x, y), (x', y'))=
\omega(x, y')-\omega(y, x')
\] 
where $(x, y, x', y')\in W$. This inner product has signature $(2n, 2n)$ and it indeed satisfies (\ref{compat}) as
\[
g((x, y), N(x', y'))=g((0, x'), (x, y))=\omega(x', x)=- g(N(x, y), (x', y')).
\]
\subsubsection{Example}
\label{example2}
Let $V$ be a $4n$--dimensional vector space over $\F=\R$ with two null structures $N_1, N_2$, such that
\[
N_1N_2+N_2N_1=-\mbox{Id},
\]
where $\mbox{Id}$ is the identity endomorphism on $V$. Then the endomorphisms
\[
I:=N_1+N_2,\qquad S:=N_1-N_2,\qquad T:=[N_1, N_2]
\]
equip $V$ with a pseudo--quaternionic structure. Indeed
\begin{eqnarray*}
I^2&=&N_1N_2+N_2N_1=-\id,\\
S^2&=&-N_1N_2-N_2N_1=\id,\\
T^2&=&N_2(-\id-N_2N_1)N_1+N_1(-\id-N_1N_2)N_2=\id,\\
IS&=&-T=-SI, \qquad IT=S=-TI, \qquad ST=I=-TS.
\end{eqnarray*}
If we instead consider $V$ over $\F=\C$ then
$I, J:=iS, K:=-iT$ form a complexified hyper--complex structure on $V$. Then $\frac{1}{2}(I\pm iJ)$ are null
structures on $V$.
\section{Null K\"ahler structures}
\label{section3}
Let $\X$ be a smooth manifold of real dimension $4n$. We shall equip $\X$ with a null  structure by smoothly extending 
such structure from each tangent space, and imposing integrability conditions.
\begin{defi}
A null structure $N$ on $\X$ is an endomorphism $N: T\X\rightarrow T\X$ such that $N^2=0$, and the sub--bundle
$\D\subset T\X$ consisting of vectors fields annihilated by $N$ has rank $2n$, and is Frobenius--integrable, i. e.
$
[\D, \D]\subset \D.
$
\end{defi}\noindent
The integrability condition holds if $N[NX, NY]=0$ or equivalently if 
$
T(NX, Y)=0
$
for all vector fields $X,Y$, where
\[
T(X, Y):=[NX, NY]-N[NX, Y]-N[X, NY]
\]
is the Nijenhuis tensor of $N$.

Null--structures are also called almost--tangent 
structures 
\cite{yano1}, and the following example shows why
\subsubsection{Example}
\label{example3}
Let $\X=TM$ be the total space of the tangent bundle to a $2n$ dimensional 
manifold $M$. Let $U\in \Gamma(TM)$ be a vector field
on $M$. Recall \cite{yano2} that the vertical lift $U^V$ of $U$ to $TM$ is a section of $T(TM)$ defined by
\[
U^V(f)=\frac{d}{d\epsilon}\bigg\rvert_{\epsilon=0} f(m, u+\epsilon U)
\]
where $f:TM\rightarrow \R$ is an arbitrary function, $m\in M$ and $u\in T_m M$. The canonical null structure on $TM$ is the endomorphism 
$N:T(TM)\rightarrow T(TM)$ defined by
\be
\label{eq5}
N(X)=[\pi_*(X)]^V,
\ee
where $\pi_*$ is the tangent map to the bundle projection $\pi:TM\rightarrow M$.

Let $(x^1, \dots, x^{2n})$ be local coordinates on $M$ covering a neighbourhood
of $m\in M$, and $(y^1, \dots, y^{2n})$ be the natural coordinates
on $T_mM$ obtained by writing a tangent--vector as $U=\sum_{i} y^i\frac{\p}{\p x^i}$. Then $(x^i, y^i)$ are local coordinates on $TM$. If
\[
X=\sum_{i} A^i\frac{\p}{\p x^i}+B^i\frac{\p}{\p y^i}
\]
is a general vector field on $TM$, then (\ref{eq5}) implies
\[
N(X)=\sum_i A^i\frac{\p}{\p y^i}.
\]
Thus, in the natural coordinate system $(x^i, y^i)$ on $TM$ the null structure is given by a $(1, 1)$ tensor\footnote{There is another, equivalent definition of this canonical null 
structure which we shall now describe. 
We have defined vertical lifts of vector fields to the tangent bundle. 
We can also define vertical lifts of functions:
if $h:M\rightarrow \R$, then $h^V=h\circ\pi$ is a function on $TM$. 
Vertical lifts of all tensor fields are  defined by
$(P\otimes Q)^V=P^V\otimes Q^V$. In particular the vertical lift of the 
identity endomorphism of $TM$ to $T(TM)$ is given by $N$.}
\[
N=\sum_{i} dx^i\otimes\frac{\p}{\p y^i}.
\]
\begin{defi}
\label{defi3}
A signature $(2n, 2n)$ pseudo--Riemannian  metric $g$ on a manifold $\X$ with a null structure $N$ is called null--K\"ahler if
\be
\label{defi_Ng}
g(NX, Y)=-g(X, NY)\qquad\mbox{and}\qquad \nabla N=0,
\ee
where $\nabla$ is the Levi--Civita connection of $g$.
\end{defi}\noindent
The fundamental two--form $\Omega\in \Lambda^2(T^*\X)$  defined by
\[
\Omega(X, Y)=g(NX, Y)
\]
is covariantly--constant, and therefore closed.
It satisfies
\[
\Omega^{\wedge n}:=\underbrace{\Omega\wedge\dots\wedge\Omega}_n\neq 0, \quad \Omega^{\wedge (n+1)}=0.
\]
This should be contrasted with the K\"ahler condition, where $\Omega^{\wedge 2n}\neq 0$, and justifies the terminology.
\subsubsection{Example}
\label{example4}
Let $\X_\C$ be a complexified hyper--K\"ahler manifold, i.e. a complex manifold of complex dimension $4n$ with
three holomorphic parallel complex structures $I, J, K$ satisfying the algebra of quaternions and Hermitian
with respect to a holomorphic metric $g$ on $\X_\C$. Then $N=\frac{1}{2}(I+iJ)$ is a (one of many) 
null--K\"ahler structure on $\X_\C$.
This example underlies the occurrence   of null structures in the geometric approach to 
Donaldson--Thomas invariants \cite{Bridgeland1, Bridgeland3}.
\vskip5pt
In the \S\ref{section31} we shall present a canonical normal--form of null--K\"ahler metrics. In the rest of this section we list properties of such metrics which do not refer to any choices of coordinates. 
\begin{prop}
The Riemann curvature $R$, the Ricci curvature $r$ and the Ricci scalar $S$ of a null--K\"ahler metric satisfy
\begin{subequations}
\label{curvatures}
\begin{align}
R(X, Y)\circ N&=N\circ R(X,Y), \label{cu1}\\
R(NX, Y)&=-R(X, NY), \label{cu2} \\
r(NX, Y)&=0, \label{cu3} \\
S&=0 \label{cu4}.		
\end{align}
\end{subequations}
\end{prop}
\noindent
{\bf Proof.}
Formula (\ref{cu1}) follows directly from $\nabla N=0$, and the definition of the curvature
\[
R(X, Y)V=[\nabla_X, \nabla_Y]V-\nabla_{[X, Y]}V.
\]
To prove (\ref{cu2}) we use the (\ref{cu1}) together with the 
symmetry properties of the Riemannian curvature:
\begin{eqnarray*}
g(R(NX, Y)V, U)&=&g(R(U, V)Y, NX)=
-g(N R(U, V)Y, X)=-g(R(U, V)NY, X)\\
&=&-g(R(X, NY)V, U).
\end{eqnarray*}
 From its definition
\[
r(X, Y)=\mbox{Tr}(V\rightarrow R(V, X)Y).
\]
The third formula (\ref{cu3}) then follows if we take $V\in \D$ as, setting 
$V=NU$, we have
\begin{eqnarray*}
r(NX, Y)&=&\mbox{Tr}(NU\rightarrow R(NU, NX)Y)\\
&=&-\mbox{Tr}(NU\rightarrow R(U, N^2X)Y)=0.
\end{eqnarray*}
Finally, to prove (\ref{cu4}) we shall compute $S=\mbox{Tr}_g(r)$ in the basis adapted to $N$,
and regard $r$ and $g^{-1}$ as linear maps.
The property (\ref{cu3}) implies that in this basis the matrix of $r$ is of the form
\[
\left(\begin{array}{cc}
*&0\\
0&0
\end{array}
\right).\]
The property (\ref{compat}) implies that the matrix of $g^{-1}$ is of the form
\[
\left(\begin{array}{cc}
0&\omega\\
\omega^T&\Theta
\end{array}
\right).
\]
for some  block $2n$ by $2n$ matrices $\Theta$ and $\omega$ such that $\omega$ is skew and
non--degenerate. Therefore $S=\mbox{Tr}(g^{-1}\cdot r)=0$.
\koniec

The next result relates null--K\"ahler structures to special holonomy, and manifolds with parallel 
spinors \cite{bryant, kath}

\begin{prop} 
\label{prop_spinor}
A null--K\"ahler manifold admits a canonical parallel pure spinor.
\end{prop}\noindent
{\bf Proof.}  This is really a result in linear algebra which builds on a bijection between the set of  pure semi-spinors in $V=\R^{2n, 2n}$ and the Grasmannian of totally
null $2n$-dimensional subspaces of $V$.
Let ${\mathcal{C}}l(2n, 2n)$ be a real $2^{4n}$--dimensional Clifford
algebra (see, e. g. \cite{ML}) of a null-Hermitian space $(V=\R^{2n, 2n}, N)$.
This algebra is generated by $2^{2n}\times 2^{2n}$ matrices
$\gamma(X)$ subject to the relations
\[
\gamma(X)\gamma(Y)+\gamma(Y)\gamma(X)=2g(X, Y){\bf 1}, \quad\mbox{where}\quad
X, Y\in V.
\]
The multiplicative group $\mbox{Spin}(2n, 2n)$ 
is generated by all elements
$\gamma(X)\gamma(Y)$, where $X, Y$ are vectors of squared norm $\pm 1$.
The {\em spin space} $\spp$ is a reducible representation space
of $\mbox{Spin}(2n, 2n)$. It can be decomposed as
\[
\spp=\spp_+\oplus\spp_-\cong\R^{2^{2n-1}}\oplus\R^{2^{2n-1}},
\]
where $\spp_{\pm}$ are irreducible spaces of {\em semi-spinors}.
It is a simple algebraic fact that any totally isotropic
subspace of $V$ has dimension at most $2n$. A semi-spinor $\iota$ is called
{\em pure} iff
\be
\label{pures}
\dim\{X\in V, \gamma(X)\iota=0\}=2n.
\ee
The  system of equations underlying (\ref{pures}) is of rank $2n$, and defines
$2n$-dimensional
plane. This plane is totally isotropic as
\[
0=\gamma(X)\gamma(X)\iota=g(X,X)\iota
\]
so $g(X, X)=0$. The space of totally isotropic planes in $\R^{2n, 2n}$ has
two components defined by a pure element of $\spp_+$ and $\spp_-$
respectively. A pure semi-spinor $\iota$ is annihilated by  
$\gamma(X_1)\gamma(X_2)\cdots\gamma(X_{2n})\in {\mathcal Cl}(2n, 2n)$ where $X_1, \cdots,
X_{2n}$ span a totally null plane. Therefore $\iota$ corresponds to an
element of the Grassmann algebra
$\xi\in \Lambda^{2n}(V^*)$ such that $\xi\wedge\xi=0$,
and the assertion of the Proposition follows because
$\xi$ is defined by a null-Hermitian structure
$
\xi=\Omega^{\wedge n},
$
and therefore is parallel.
\koniec
If $n=1$ then all semi-spinors are pure. The first non-trivial case
corresponds to 8-dimensional null-Hermitian structures.
\subsection{Null K\"ahler potential}
\label{section31}
In this section we shall find a canonical normal form of a null--K\"ahler metric, and
express it in terms of second derivatives of one arbitrary function on $\X$.

\noindent
{\bf Proof of Theorem \ref{theorem1}.} Let
\[
M=\X/\mbox{Ker}(N)
\]
be a $2n$--dimensional quotient manifold by the kernel of the $2n$--dimensional integrable distribution $\D$
of vector fields annihilated by $N$. Locally we regard $\X$ is the tangent bundle to $M$,
and in the coordinate system of Example \ref{example3} the endomorphism $N$ is
\be
\label{latestN}
N=\sum_{i} dx^i\otimes\frac{\p}{\p y^i},
\ee
where $(x^1, \dots, x^{2n})$ are local coordinates on $M$, and $(y^1, \dots, y^{2n})$ are linear
coordinates on fibres of $TM\rightarrow M$. The distribution 
$\D=\mbox{span}\{\p/\p y^1, \dots, \p/\p y^{2n}\}$ is totally isotropic, and therefore
there exists a collection of functions on $\X$
\[
\omega_{ij}=\omega_{ij}(x, y), \quad \Theta_{ij}=\Theta_{ij}(x, y)
\]
such that $\Theta_{ij}=\Theta_{ji}$ and
\be
\label{metric_meantime}
g=\frac{1}{2}\sum_{i, j}\omega_{ij}(dy^i\otimes dx^j+dx^j\otimes dy^i)+
\frac{1}{2}\sum_{i, j}\Theta_{ij}(dx^i\otimes dx^j+dx^j\otimes dx^i).
\ee
Evaluating the null--K\"ahler condition (\ref{defi_Ng}) on coordinate vector fields shows that 
\[
\omega_{ij}=-\omega_{ji}.
\]
We now impose the parallel condition $\nabla N=0$, where $N$ is given by (\ref{latestN}). The
$dy^i\otimes dx^j$ components of $\nabla N$ vanish iff
\[
\frac{\p \omega_{ij}}{\p y^k}=0
\]
so that $\omega_{ij}=\omega_{ij}(x)$, and
\[
\omega=\frac{1}{2}\sum_{i, j}\omega_{ij} dx^i\wedge dx^j
\]
is a symplectic form on $M$. Locally there exists a diffeomorphism of $M$ to a Darboux coordinate system $\tilde{x}=\tilde{x}(x)$ such 
that
\[
\omega=d\t{x}^1\wedge d\t{x}^{n+1}+\dots+ d\t{x}^n\wedge d\t{x}^{2n}.
\]
The transformation of $\X=TM$ induced by this diffeomorphism is
\[
\t{x}^i=\t{x}^i(x^j),\quad \t{y}^i=\sum_{j}\frac{\p\t{x}^i}{\p x^j}y^j, \quad i, j=1, \dots, 2n,
\]
and it 
preserves the form of $N$, as $N=\sum_i d\t{x}^i\otimes\p/\p \t{y}^i$. We shall use this Darboux coordinate system from now on, and drop tildes.
This has an effect of reducing $\omega_{ij}$ to a constant symplectic matrix as in (\ref{nk_form}).
We now move on the vanishing of $dx^i\otimes dx^j$ 
components in $\nabla N$. This is equivalent to
\[
\frac{\p \Theta_{ij}}{\p y^k}=\frac{\p \Theta_{ik}}{\p y^j}
\]
which gives the integrability conditions for the existence of $2n$ function $\Theta_1, \dots, \Theta_{2n}$ on $\X$ such that
\[
\Theta_{ij}=\frac{\p\Theta_i}{\p x^j}.
\]
The symmetry condition $\Theta_{ij}=\Theta_{ji}$ implies the existence of a single function $\Theta=\Theta(x, y)$ such that
\[
\Theta_i=\frac{\p\Theta}{\p y^i}.
\]
This puts the metric (\ref{metric_meantime}) in the canonical form (\ref{nk_form}).
\koniec
The local normal form (\ref{nk_form}) is not invariant under general diffeomorphisms of $\X$. The subgroup
of the pseudogroup of all diffeomorphisms changing coordinates $(x^i, y^i)$ as well as $\Theta$, while
preseving (\ref{nk_form}) is a semi--direct product of $\mbox{SDiff(M)}$ and $\Gamma(M)$, where the symplectomorphisms 
$\mbox{SDiff}(M)$ of $M$ act on $\X=TM$ by a Lie lift, and $\Gamma(M)$ acts on the fibres of $TM$ by translations. The details are as follows:
Let $Y$ be a vector field on $\X$ generating a one--parameter group of diffeomorphisms. The conditions
\[
{\mathcal L}_Y\Omega=0, \quad {\mathcal L}_YN=0
\] 
imply
\[
Y=\sum_{i, j} \omega^{ij}\frac{\p H}{\p x^j}\frac{\p}{\p x^i}+
\sum_{i, j, k}
\Big(y^k\omega^{ij}\frac{\p^2 H}{\p x^k \p x^j}\Big) 
+\sum_{i}T^i\frac{\p}{\p y^i},
\]
where $\omega^{ij}$ is the inverse matrix of $\omega_{ij}$, 
i. e.
$\omega^{ik}\omega_{kj}={\delta^i}_j$,
and
$(H, T^1, \dots, T^{2n})$ are arbitrary functions of $(x^1, \dots, x^{2n})$. 
Set
\[
\tilde{x}^i=x^i+\epsilon {\mathcal L}_Y(x^i) ,\quad \tilde{y}^i=y^i+\epsilon {\mathcal L}_Y(y^i)
, \quad \widetilde{\Theta}=\Theta+\epsilon\delta\Theta.
\]
Using
\[
\frac{\p}{\p \tilde{y}^i}=\frac{\p}{\p y^i}-\sum_{j, k}\epsilon\omega^{jk}\frac{\p^2 H}{\p y^i\p y^k}\frac{\p}{\p y^j}
\]
we find that the action generated by $Y$ preserves the form of $g$, i. e.
\[
{g}=
\frac{1}{2}\sum_{i, j} \omega_{ij}(d\tilde{y}^i\otimes d\tilde{x}^j+d\tilde{x}^j\otimes d\tilde{y}^i) 
+\frac{\p^2 \widetilde{\Theta}}{\p \tilde{y}^i \p \tilde{y}^j}(d\tilde{x}^i 
\otimes d\tilde{x}^j +d\tilde{x}^j \otimes d\tilde{x}^i)+O(\epsilon^2)
\]
if 
\be
\label{freedom}
\delta \Theta=\sum_{i, j, k} \Big(\frac{1}{6}y^iy^jy^k\frac{\p^3 H}{\p x^i\p x^j \p x^k}-\frac{1}{2}y^jy^k\omega_{ij}
\frac{\p T^i}{\p x^k}\Big)+\sum_i y^i Q_i+R,
\ee
where $(Q_1, \dots, Q_{2n}, R)$ are arbitrary functions of $x^i$.
\subsection{Null--K\"ahler Einstein metrics}
\label{section32}
Computing the Ricci tensor of a null--K\"ahler metric in the form (\ref{nk_form})  we find
\[
r=\sum_{i, j}\frac{\p^2 f}{\p y^i\p y^j} (dx^i\otimes dx^j+dx^j\otimes dx^i),
\]
where
\be
\label{form_of_f}
f\equiv\sum_{i, j}\omega^{ij}\frac{\p^2 \Theta}{\p y^i\p x^j}+\sum_{i, j, k, l}\frac{1}{2}\omega^{ik}\omega^{jl}
\frac{\p^2 \Theta}{\p y^i\p y^j} \frac{\p^2 \Theta}{\p y^k\p y^l}.
\ee
The Ricci--flat condition on $g$ therefore reduces to a system of fourth order PDEs on $\Theta$
which can be integrated twice to give a single second orde PDE on $\Theta$
\be
\label{einstein_eq}
f=G+y^iF_i,
\ee
where 
$(G, F_1, \dots, F_{2n})$ are arbitrary functions of $x^i$. Applying the Cauchy--Kovalevskaya theorem shows that in the real--analytic category
the general Ricci--flat null K\"ahler metric depends on two arbitrary functions of $4n-1$ 
variables and some number of functions of $2n$ variables.
\subsubsection{Example} It can be explicitly verified
that for
\be
\Theta=\frac{c}{\rho^{2n-1}}\quad\mbox{where}\quad  \rho= \sum_{i, j} \omega_{ij}y^i x^j, \quad
c=\mbox{const}
\ee
the linear and non--linear terms in (\ref{form_of_f}) vanish separately resulting in $f=0$. The resulting metric
\[
g=\frac{1}{2}\sum_{i, j} \omega_{ij}(dy^i\otimes dx^j+dx^j\otimes dy^i)
+\frac{2c n(2n-1)}{\rho^{2n+1}}\Big(\sum_{k, l}\omega_{kl}
x^k dx^l\Big)^{\otimes 2}
\]
is therefore null--K\"ahler, and Ricci--flat. This metric with $n=1$ is the 
Sparling--Tod $H$--space \cite{ST}.
\subsection{Complex hyper--K\"ahler metrics with affine symplectic fibrations}
\label{section33}
In the complexified setting the coordinates $(x^i, y^i)$ are holomorphic on
the complex manifold $\X_\C$ of complex dimension $4n$. If
$\Theta=\Theta(x, y)$ satisfies the system of PDEs
\begin{eqnarray}
\label{BS}
H_{ij}&=&0,\quad i, j=1, \dots, 2n \quad\mbox{where}\nonumber\\
H_{ij}&\equiv&\frac{\p^2 \Theta}{\p y^i \p x^j}-\frac{\p^2 \Theta}{\p y^j \p x^i}+
\sum_{k, l}\omega^{kl}\frac{\p^2 \Theta}{\p y^i \p y^l} \frac{\p^2 \Theta}{\p y^j \p y^k}=0,
\end{eqnarray}
then the metric (\ref{nk_form}) is complexified hyper--K\"ahler. In 
\cite{Bridgeland3} it was shown that if $(M_\C, \omega)$ is a complexifed
symplectic manifold of complex dimension $2n$, and $g$ is a complexified 
hyper--K\"ahler metric on $\X_\C=TM_\C$ such that the null--K\"ahler two--form 
\[
\pi^*(\omega)=\Omega_I+i\Omega_J,
\] 
then $g$ is locally
of the form (\ref{nk_form}), where $\Theta$ satisfies the system (\ref{BS}).

The system (\ref{BS}) consists of some of the flows of the hyper--K\"ahler hierarchy
\cite{Tak, DM}. It implies the Frobenius integrability
\be
\label{thels}
[l_i, l_j]=0, \quad i, j=1,\dots, 2n
\ee
for the rank--$2n$ distribution 
spanned by
\[
l_i=\frac{\p}{\p y^i}+\lambda\Big(\frac{\p}{\p x^i}
+\sum_{j, k}\omega^{jk}\frac{\p^2 \Theta}{\p y^i\p y^j}\frac{\p}{\p y^k}
\Big)
\]
on $\X_\C\times \CP^1$, where $\lambda$ is the affine coordinate on $\CP^1$.

Vanishing of the Lie brackets (\ref{thels}) gives a weaker set of 
conditions\footnote{Note that $\sum_{i, j}\omega^{ij}H_{ij}=2f$. Therefore 
(\ref{weaker_form}) implies Ricci--flatness, but the converse is not true.}
\be
\label{weaker_form}
\frac{\p H_{ij}}{\p y^k}=0, \quad \mbox{so that} \quad H_{ij}=C_{ij}(x),
\ee
for some skew $C_{ij}$. If $n=1$ then a transformation $\Theta\rightarrow 
\Theta+\sum_{i}y^iQ_i(x)$
can be used to set $C_{ij}$ to zero. For general $n$ this only seems possible
if the two--form $\sum_{i, j}C_{ij}dx^i\wedge dx^j$ is closed. In \cite{Bridgeland3}
it has been argued that for any $n$ the conditions (\ref{weaker_form})
together with the additional assumption that the function $\Theta$
is odd in the fibre variables $(y^1, \dots, y^{2n})$ imply (\ref{BS}).

 Geometrically, $\lambda$ labels the $2n$--dimensional surfaces (the $\alpha$--surfaces
in the twistor approach \cite{DM}) through each point of $\X_\C$. The twistor 
space of $(\X_\C, g)$ is the space of these $\alpha$--surfaces. It is
a complex manifold of complex dimension $2n+1$, which
arises as the quotient of $\X_\C\times \CP^1$ by the distribution spanned by 
$l_i$.
The points in $\X_\C$ correspond to rational curves in 
${\mathcal Y}$ with normal bundle $\C^{2n}\otimes \OO(1)$, where
$\OO(1)$ is the line bundle with Chern class 1 on $\CP^1$. In \cite{DM}
this twistor correspondence has been extended to the full hyper--K\"ahler 
hierarchy. 
\subsubsection{Example} A strong Joyce structure has, in \cite{Bridgeland2},
been defined to be a solution $\Theta$ of the system (\ref{BS}) subject to three additional conditions:
\begin{enumerate}
\item $\Theta$ is odd in the variables $y^i$.
\item $Z\equiv \sum_ix^i\frac{\p}{\p x^i}$ is a homothetic Killing vector field
such that
\[
{\mathcal L}_Z g=g, \quad {\mathcal L}_Z \Theta=-\Theta.
\]
\item  The metric is invariant under the lattice transformations
\[
y^i\rightarrow y^i+2\pi\sqrt{-1}, \quad i=1, \dots, 2n.
\]
\end{enumerate}
An example of a solution to (\ref{BS}) which also satisfies the three conditions above is
\be
\Theta=\frac{\sinh{y^1}}{x^1}.
\ee
The resulting metric is non--flat, and is an example of a Ricci--flat plane wave.
\subsection{Conformal invariance}
In four dimension the restricted conformal transformations, where the 
conformal factor is constant along the distribution 
${\mathcal D}=\mbox{Ker}(N)$,
preserve the null--K\"ahler condition: If $F=F(x^1, x^2)$, and
\[
\hat{g}=F^2g, \quad\hat{\Omega}=F^3\Omega, \quad\mbox{then}\quad
\hat{\nabla}\hat{\Omega}=0.
\] 
This conformal invariance is not present in other dimensions: for
$F^k\Omega$ to be closed we need $k=0$, and then $\hat{\nabla}\Omega=0$ implies
$F=\mbox{const}$.
\subsection{Walker structures}
Recall that a distribution ${\mathcal D}$ on a pseudo--Riemannian manifold 
$\X$ is called parallel if $\nabla_Y X\in \Gamma({\mathcal D})$ for all
$X\in \Gamma({\mathcal D})$ and $Y\in \Gamma(T\X)$. The pseudo--Riemannian manifolds admitting a parallel distribution of rank equal to half of the manifold 
dimension are called the Walker manifolds \cite{Gilkey}, and it was shown by Walker
\cite{Walker_1}, that locally a Walker metric is of the form (\ref{metric_meantime}) for some
functions $\Theta_{ij}$. The null--K\"ahler manifolds form a subclass of the Walker manifolds,
where $\Theta_{ij}$ is a Hessian of one function. Indeed, any vector field in ${\mathcal D}$
is of the form $N(X)$ for some $X\in \Gamma(T\X)$, and we have
\[
\nabla_Y N(X)=N \nabla_Y X\in\mbox{Ker}(N).
\]
\section{Four dimensions}
\label{section4}
The author has  first came across the null--K\"ahler structures when investigating twistor theory and integrability of a certain fourth order PDE in
four dimensions \cite{D02}. 
In four dimensions the existence of  a maximal rank parallel endomorphism $N$ with $N^2=0$ is equivalent to existence of a parallel semi--spinor, 
i. e. a parallel section of a rank--two symplectic vector bundle which we chose to be $\spp_+\rightarrow \X$ where
\be
\label{can_bun_iso}
T \X\cong {\spp_+}\otimes {\spp_-}
\ee
is a canonical bundle isomorphism, and $\spp_-$ is another rank--2 symplectic vector bundle. 
This isomorphism is related to the metric on $\X$ by
\[
g(v_1\otimes w_1,v_2\otimes w_2)
=\varepsilon_+(v_1,v_2)\varepsilon_-(w_1, w_2)
\]
where $v_1, v_2\in \Gamma(\spp_+)$ and $w_1, w_2\in \Gamma(\spp_-)$, and $\varepsilon_{\pm}$ are symplectic  
structures on $\spp_\pm$ which are parallel with respect to $\nabla$.
The Hodge $\ast$ operator is an
involution on two-forms, and induces a decomposition
\be \label{splitting}
\Lambda^{2}(T^*\X) = \Lambda_{+}^{2}(T^*\X) \oplus \Lambda_{-}^{2}(T^*\X)
\ee
of two-forms
into self-dual (SD)
and anti-self-dual (ASD)  components. 
 Given a parallel section of $\spp_+$, another isomorphism
\be
\label{another_iso}
\Lambda_{+}^{2}\cong \mbox{Sym}^2({\spp_+}^{*})
\ee
implies the existence of a parallel self--dual two--form $\Omega$ such that $\Omega\wedge\Omega=0$.
This two--form, together with $g$ define the nilpotent endomorphism $N$ by (\ref{hermitean_form}).

In four dimensions there are three non--linear systems of PDEs, one of them completely solvable, one integrable, and one not-integrable, which can be imposed on 
the null structure. Before writing these systems down in coordinates of Theorem (\ref{theorem1}) recall \cite{AHS}
that in four dimensions the Riemann tensor of $g$ can be regarded  
as a map $\mathcal{R}: \Lambda^{2}(T^*\X) \rightarrow \Lambda^{2}(T^*\X)$
which admits a decomposition   under the splitting  (\ref{splitting}):
\be \label{decomp}
{\mathcal R}=
\left(
\mbox{
\begin{tabular}{c|c}
&\\
$C_+ +\frac{1}{12}S$&$r_0$\\ &\\
\cline{1-2}&\\
$r_0$ & $C_-+\frac{1}{12}S$\\&\\
\end{tabular}
} \right) .
\ee
Here $C_{\pm}$ are the SD and ASD parts 
of the Weyl tensor, $r_0$ is the
trace-free Ricci curvature, and $S$ is the scalar curvature which acts
by scalar multiplication.
We are now ready to present the three systems of PDEs 
\subsection{Self--dual null--K\"ahler} The condition $C_-=0$ is equivalent to
\[
\frac{\p^4\Theta}{\p y^i\p y^j\p y^k \p y^l}=0.
\]
Therefore the most general self--dual null--K\"ahler metric in four dimensions is of the form
(\ref{nk_form}) with 
\[
\Theta=\sum_{i, j, k}\Gamma_{ijk}y^iy^jy^k
\]
where the functions $\Gamma_{ijk}$ depend only on $(x^1, x^2)$, as the coordinate freedom (\ref{freedom}) can be used to remove the quadratic
and linear terms from $\Theta$. The resulting metric is a Walker's projective extension 
\cite{DT_k} of a projective structures on the surface $M=\X/\D$.
\subsection{ Anti--self--dual null--K\"ahler} The condition $C_+=0$ is equivalent to
a 4th order PDE for $\Theta$:
\begin{eqnarray}
\label{eqmd}
f&=& \Theta_{x^1y^2}-\Theta_{x^2y^1}+\Theta_{y^1y^1}\Theta_{y^2y^2}-{(\Theta_{y^1y^2})}^2\\
\Delta_g f&:=&f_{x^1y^2}-f_{x^2y^1}+\Theta_{y^2y^2}f_{y^1y^1}+\Theta_{y^1y^1}f_{y^2y^2}
-2\Theta_{y^1y^2}f_{y^1y^2}=0\nonumber
\end{eqnarray}
where  subscripts denote partial derivatives, i. e.  $\Theta_{x^1}=\p\Theta/\p x^1$ etc.
Note that $\Delta_g$ is the Laplace--Beltrami operator of the metric $g$, and the expression for $f$ agrees with the general formula (\ref{form_of_f}). This equation is integrable by twistor transform \cite{D02}, the dressing method
\cite{Bogdanov}, and the Manakov--Santini inverse scattering transform \cite{Yi}. 
The general solution depends on
4 functions of 3 variables.
It has recently been shown to
arise from a second-order integrable Lagrangian \cite{Ferapontov}. 
\subsection{Null--K\"ahler Einstein}
As the scalar curvature of null--K\"ahler manifolds always vanishes,
the Einstein condition is equivalent to the vanishing of the Ricci tensor of
 $g$. The resulting second order PDE (\ref{einstein_eq}) on $\Theta$ is the hyper–heavenly equation of
Pleba\'nski and Robinson \cite{pleban_robinson} for non--expanding metrics with self--dual Weyl tensor
$C_+$ of type ${\bf N}$. (Recall that $(\X, g)$
is called hyper–heavenly if the self-dual Weyl tensor is algebraically special, i. e. has a repeated root when
regarded as a binary quartic. Type ${\bf N}$ corresponds to a repeated root of order $4$).
\subsection{Heavenly equation} Imposing the Einstein condition together with the anti--self--duality of
the Weyl tensor reduces the 4th order equation (\ref{eqmd}) to a second order PDE
\[
f=0.
\]
This is Pleba\'nski's  second heavenly equation \cite{pleban}. The resulting
metric is pseudo--hyper--K\"ahler.
\section{Anti--self--duality and isomonodromy}
\label{section5}
In this section we shall assume that 
$(\X, N, g)$ is a null--K\"ahler four--manifold with anti--self--dual Weyl curvature which is cohomogeneity-one, i. e. there exists an isometry group $G$ acting transitively on three--dimensional surfaces in $\X$. 
The four-dimensional cohomogeneity-one metrics can be classified according to the Bianchi type of the three-dimensional real Lie algebra of $G$. Locally $\X=\R\times G$ and the ASD cohomogeneity--one
null-K\"ahler condition reduces to solving a system of ODEs. To write this system down, and recognise it as
the isomonodromy problem for Painlev\'e I and II if $G=SL(2)$, we shall use the twistor methods 
\cite{penrose_NG, D02}.
We shall therefore work in the holomorphic category, and assume that $\X_\C$ is a complex oriented four--manifold, and $(N, g)$ are holomorphic.
\begin{defi}
\label{alpha_defi}  
An $\alpha$-surface is  two-dimensional
surface $\zeta\subset \X_\C$ such that for all $p \in \X_\C$ the tangent space 
$T_p\xi$ is a totally null plane with self--dual tangent bi-vector.
\end{defi}

The Nonlinear Graviton theorem of Penrose \cite{penrose_NG} states  that
there locally exist a three--parameter family of $\alpha$--surfaces
iff the self-dual part of the Weyl tensor of $g$ vanishes. The twistor space ${\mathcal Y}$
of an ASD four-manifold is defined to be the space of $\alpha$--surfaces. 
It is a three--dimensional complex manifold with a four--parameter family of
rational curves with normal bundle $\OO(1)\oplus\OO(1)$. The points of the three--dimensional twistor
space ${\mathcal Y}$ are $\alpha$--surfaces in $\X_\C$: There is a rational curve
$L_p\cong\CP^1$ worth of such surfaces  through each point $p\in \X_\C$, and therefore 
points in $\X_\C$ correspond to rational curves in ${\mathcal Y}$. 
The conformal structure on $\X_\C$ is defined by declaring two points $p_1, p_2 \in \X_\C$ to be null--separated
iff the corresponding rational curves in ${\mathcal Y}$ intersect at one point.
 
The correspondence between $\X_\C$ and ${\mathcal Y}$
can be expressed in the  the double fibration picture (see e.g. \cite{mason}).
\be
\label{doublefib}
{\X_\C}\stackrel{r}\longleftarrow 
{\mathcal F}\stackrel{q}\longrightarrow {\mathcal Y},
\ee
where the five--complex--dimensional 
correspondence is defined by
\[
{\mathcal F}={\mathcal{Y}}\times {\X_\C}|_{\zeta\in L_p}= {\X_\C}\times\CP^1
\] 
where $L_p$ is the rational curve in ${\mathcal Y}$ that corresponds to 
$p\in {X_\C}$, and $\zeta\in{\mathcal{Y}}$ lies on $L_p$. 
The twistor space arises
as a quotient of ${\mathcal F}$ by a two--dimensional integrable distribution spanned by the vector fields
\be 
\label{tetradlax} 
l_1 = E_{11} -  \lambda E_{12} + f_1 \frac{\p}{\p \lambda}, 
\quad l_2 = E_{21} -  \lambda E_{22} + f_2 \frac{\p}{\p \lambda}, 
\ee
where $\lambda$ is an affine coordinate on $\CP^1$,
the  functions  $f_1, f_2$ on ${\mathcal F}$ are cubic in $\lambda$,
and  $E_{ij}$ are four independent 
holomorphic vector fields on $\X_\C$ such that the conformal structure defined by 
the contravariant metric
\be
\label{tetrad_0}
g=\frac{1}{2}(E_{11}\otimes E_{22} +E_{22}\otimes E_{11} - E_{12}\otimes E_{21}-E_{21}\otimes E_{12}).
\ee
The Frobenius  integrability condition
\be
\label{intlm}
[l_1, l_2]=0 \quad (\mbox{mod}\; l_1, l_2)
\ee
is equivalent to the anti--self--duality condition $C_+=0$ on $\X_\C$.
If the integrability condition holds then there is a $\CP^1$--worth of 
$\alpha$-surfaces spanned by $\{ E_{11} -  \lambda E_{12}, 
E_{21} -  \lambda E_{22}\}$ through any point in $\X_\C$. If all vectors $E_{ij}$ are real then the signature of $g$ is $(2, 2)$, and there exists  an $\RP^1$--worth of real $\alpha$-surfaces through each point of a real four--manifold $\X$.
The null--K\"ahler condition on top of anti--self--duality gives rise to an additional structure on the twistor space:
\begin{theo}\cite{D02}
\label{mytheo}
Let ${\mathcal Y}$ be a three-dimensional complex manifold with
\begin{enumerate}
\item A four-parameter family of rational curves with normal bundle $\OO(1)\oplus \OO(1)$.
\item  A preferred section of $\kappa^{-1/4}$ where $\kappa$ 
is the holomorphic canonical bundle of ${\mathcal Y}$.
\item An anti-holomorphic involution $\rho: {\mathcal Y}\rightarrow
{\mathcal Y}$ fixing a real equator of each
rational curve, and leaving the section of $\kappa$ above invariant.
\end{enumerate}
Then the real moduli space $\X$ of the $\rho$--invariant curves is equipped with a restricted
conformal class $[g]$ of anti--self--dual null-K\"ahler metrics: 
if $g\in [g]$, and $\Omega$ is a null-K\"ahler
two-form, then $(\hat{g}=F^2g, \hat{\Omega}=F^3\Omega)$ is also null--K\"ahler
 for any function $F$ such that $dF\wedge\Omega=0$.

 Conversely, given
a real analytic ASD null--K\"ahler metric, there exists a corresponding twistor space
${\mathcal Y}$
with the above structures.
\end{theo}
If one  is only interested in the complexified picture, where $g$ and $N$ are
holomorphic on $\X_\C$, then condition (3) in Theorem \ref{mytheo} can be dropped.

From now on we shall additionally assume that there exists
a three--dimensional complex Lie group $G$ acting on $\X_\C=\C\times G$ by 
isometries with generically three--dimensional orbits.
We shall make a choice for $G$, and take it to be $SL(2, \C)$ (or $SL(2, \R)$ if $\X$ is a real four--manifold
with a $(2, 2)$ metric).  Its Lie algebra is generated by the left invariant vector fields
$L_1, L_2, L_3$ on $G$ which satisfy
\be
\label{bianchi_2}
[L_1, L_2]=L_2, \quad [L_1, L_3]=-L_3, \quad [L_2, L_3]=2L_1.
\ee
The conformal isometries are generated by 
the right-invariant vector fields
$R_\alpha, \alpha=1, 2, 3$ on $G$. The metric on $\X_\C$ will be expressed in terms of the left--invariant one--forms $\sigma^1, \sigma^2, \sigma^3$ on $SL(2)$ such that
\[
{\mathcal L}_{R_{\alpha}}\sigma^\beta=0, \quad
L_\alpha\hook\sigma^\beta=\delta^{\beta}_{\alpha}
\]
and
\be
\label{one_formss}
d\sigma^1=2\sigma^3\wedge \sigma^2, \quad d\sigma^2=\sigma^2\wedge \sigma^1,\quad
d \sigma^3=\sigma^1\wedge\sigma^3.
\ee

The $G$--action on $\X_\C$
maps $\alpha$--surfaces to $\alpha$--surfaces 
and thus gives rise to a holomorphic group action of $G$ on 
the twistor space ${\mathcal Y}$. Let the  
$\widetilde{R}_\alpha, \alpha=1, 2, 3$ be holomorphic vector
fields on ${\mathcal Y}$ generating this action and corresponding to $R_\alpha$.
Consider a quartic
\be 
\label{quartic} s = \mbox{vol}_{\mathcal Y} ({\widetilde{R}}_1,  {\widetilde{R}}_2,
{\widetilde{R}}_3), 
\ee
where  $\mbox{vol}_{\mathcal Y}$ is a holomorphic volume  form on the twistor space with values
in $\OO(4)$.
This quartic  vanishes at each twistor line at four points, where the holomorphic vector fields corresponding to the isometries become linearly dependent. 
We shall, form now on assume that
the four zeros of the quartic coincide, and so $s$ gives a preferred section of $\kappa^{-1/4}$. 
The corresponding conformal structure therefore contains a null-K\"ahler structure by 
Theorem \ref{mytheo}. We shall first need to establish two technical results about the quartic 
(\ref{quartic}), as the canonical form of the metric depends on whether $s$ vanishes identically, or not.
\begin{prop}
\label{quartic_lemma}
If the quartic (\ref{quartic}) vanishes identically then the conformal class
containing $g$ is hyper--complex, or equivalently if there exists a holomorphic fibration of ${\mathcal Y}$ over $\CP^1$ such that the twistor
curves are sections of this fibration.
\end{prop}
\noindent
{\bf Proof.}
 We shall first introduce some notation. Let $\pi^i=[\pi^1, \pi^2]$ be homogeneous coordinates on $\CP^1$--fibres of the bundle $\PP(\spp_+)$ 
such that $\lambda=-\pi^2/\pi^1$ in the patch where $\pi^1\neq 0$, and let $\pi_i=\sum_{j}\varepsilon_{ji}\pi^j$. 
Assemble the frame in (\ref{tetradlax}) into a vector--valued two by two matrix
$E$ with components $E_{ij}$
so that the twistor distribution ${\mathcal D}_{\mathcal Y}\equiv\OO(-1)\otimes\C^2$ given by
(\ref{tetradlax}) takes the form
\[
l_i=\sum_j \pi^j E_{ij}+f_i\frac{\p}{\p \lambda}, \quad i=1, 2.
\]
For this to be homogeneous of degree $1$ in $\pi$ the 
functions $(f_1, f_2)$ need to be sections
of $\OO(3)$.  Let $e^{ij}$ be a frame of one--forms dual to $E_{ij}$  so that $E_{ij}\hook e^{mn}={\delta_i}^m{\delta_j}^n$,
and the metric is given by
\[
g=\frac{1}{2}(e^{11}\otimes e^{22}+e^{22}\otimes e^{11}-e^{12}\otimes e^{21}-e^{21}\otimes e^{12}).
\]
In the double fibration picture (\ref{doublefib}) the quartic $s$  pulls back
to a quartic on ${\mathcal F}$ given by
\be 
\label{quartic1} q^*(s) = (d \lambda \wedge \mbox{vol}) (l_1, l_2,  {R}_1,  {R}_2,
{R}_3), 
\ee
where $\mbox{vol}$ is the holomorphic volume form on $\X_\C$ such that
$\mbox{vol}( E_{11}, E_{21}, E_{12}, E_{22})=1$,
and we have chosen to work in an invariant frame, where
the lifts of $R_\alpha$s to the correspondence space ${\mathcal F}$ are given by $R_\alpha$s.
Such a frame always exists, as given a cohomogeneity--one metric 
of the form
(\ref{tetrad_0}) we
can choose a frame of one--forms $e^{ij}$ which are linear combinations of the left invariant one--forms on 
on $G$, and $dt$. Here $t$ is a function
on such that the surfaces of constant $t$ are the orbits of $SL(2)$, and so $dt$is normal to the surfaces of homogeneity.
The coefficients
of this combination only depend on $t$, so the self--dual two--forms
\[
e^{11}\wedge e^{21},\quad e^{11}\wedge e^{22}-e^{21}\wedge e^{12}, \quad e^{12}\wedge e^{22} 
\]
constructed from the frame $e^{ij}$
are also $G$--invariant. Therefore the lift of the $SL(2)$ action to the bundle $\spp_+$ is trivial, but
the correspondence  space is the projectivisation  of this bundle.

 Define $T_{ij}$ by
\be
\label{formtt}
dt=\sum_{ij} T_{ij} e^{ij}, \quad  \mbox{so that}\quad T_{ij}=E_{ij}(t).
\ee
where we have used $d=\sum_{i, j} e^{ij}\otimes E_{ij}$. 
Let $\mbox{vol}=dt\wedge \mbox{vol}_{SL(2)}$, so that \[\mbox{vol}_{SL(2)}(R_1, R_2, R_3)=1.\]
Therefore,  using (\ref{formtt}),
\[
q^*(s)=\frac{1}{2}\sum_{i, j, m, n}{\varepsilon_-}^{mn}T_{ij} d\lambda\wedge e^{ij}(l_m, l_n)
\]
where $\varepsilon_-$ is the symplectic structure on ${\spp_-}^*$. Using
$e^{ij}(l_m, \cdot)={\delta^i}_m\pi^j$ gives
\begin{eqnarray}
\label{new_form}
q^*(s)&=& \sum_{i, j, k, m, n}T_{ij}f_m\pi^k{\delta^i}_n{\delta^j}_k 
{\varepsilon_-}^{mn}\nonumber\\
&=&  \sum_{i, j, k} f_i T_{jk}\pi^k{\varepsilon_-}^{ij}.
\end{eqnarray}
If the invariant frame $E_{ij}$ is also such  that $f_1=f_2=0$, then $q^*(s)$ given by (\ref{new_form}) is identically zero. In \cite{D99} 
it was shown that a frame with $f_1=f_2=0$ (and therefore a holomorphic fibration ${\mathcal Y}\rightarrow \CP^1$ \cite{Boyer, Hitchi_frame, Calderbank})
exists iff $g$ is hyper--complex. Therefore the hyper--complex condition is  necessary
for the vanishing of $s$.
\koniec
\newpage
{\bf Remarks}
\subsection{} A complexified hyper--Hermitian (which  in four dimensions is equivalent to complexified hyper--complex) structure on $\X_\C$ is a triple of holomoprhic
Hermitian endomorphisms $I, J, K$ of  $T\X_\C$ which satisfy the algebra of quaternions.
If $J=iS, K=-iT$,  and $I, S, T$ are all real, then they form a pseudo--hyper--Hermitian
 structure on split--signature real  four--manifold $\X$. The endomorphism $I$ endows $\X$ with
the structure of a two--dimensional complex K\"ahler manifold, and so does every other complex structure $aI+bS+cT$ parametrised by the points of the hyperboloid
\be 
\label{hyperboloid}
a^2-b^2-c^2=1.
\ee
\subsection{}
\label{sub51}
 The converse of Proposition \ref{quartic_lemma}
does not hold: if $g$ is ASD and  Ricci--flat (and therefore hyper--complex) but the $SL(2)$ action
rotates the covariantly constant self--dual two forms, then $s$ does not vanish. That is to say
the covariantly constant frame does not have to be invariant. The basis of two--forms in
an invariant frame (which, as we have argued, always exists) is not covariantly constant and
so $f_1$ and $f_2$ will not vanish. ASD Taub--NUT or the Atiyah--Hitchin metrics are both examples
illustrating this phenomenon. We can however say more if the isometric group action preserves the null--K\"ahler structure:
\begin{lemma}
If $(\X_\C, g, N)$ is a cohomogeneity--one $SL(2)$ invariant null--K\"ahler
structure which is Ricci flat, and such that $N$ is preserved by the group
action, then the quartic $s$ vanishes identically.
\end{lemma}
\noindent
{\bf Proof.}
Let $\iota\in \Gamma(\spp_+)$
be the covariantly constant spinor defining $N$. Then $\iota$ must be in a linear combination of the covariantly constant basis
of $\spp_+$ (which exists for ASD metrics iff they are Ricci flat)
 with constant coefficients (or it can not be parallel). Therefore the null structure $N$ belongs to the hyperboloid (\ref{hyperboloid}) of complex
structures defined by the covariantly constant basis. The group $SL(2)$ acts on this hyperboloid, and we require that
it fixes $N$. But this implies that it must fix all other points of the hyperboloid (as otherwise
the Lie algebra  relations  would  be violated). Therefore the covariantly constant  frame is also invariant
and $s=0$.
\koniec

\subsection{}Proposition \ref{quartic_lemma} was established by Hitchin who used representation--theoretic arguments 
\cite{Hitchinis} under an additional assumption that the twistor space admits a
real structure which singles out a Riemannian real section of $\X_\C$. 
In these circumstances the quartic $s$ either vanishes identically, or it admits two repeated roots, or all four roots are distinct. This assumption
is not valid in the context of null-K\"ahler structures and split signature metrics.

\begin{lemma}
If $(\X_\C, g)$ is hyper--Hermitian, and null--K\"ahler, then the metric
$g$ is conformal to a Ricci--flat metric.
\end{lemma}
\noindent
{\bf Proof.} We shall use the formulation of the hyper--Hermitian condition due to Boyer \cite{Boyer}, which is also applicable
in the complexified setting \cite{D99}: a metric on $\X_\C$ is hyper--Hermitian if and only if
there exists a basis $(\Sigma^1, \Sigma^2, \Sigma^2)$ of ${\Lambda^2}_+$ such that
\be
\label{boyer_hc}
d\Sigma^{\alpha}=2A\wedge\Sigma^{\alpha}, \quad \alpha=1, 2, 3
\ee
for some one--form $A$. Moreover a hyper--Hermitian $g$  is locally  conformal to Ricci--flat iff $A$ closed.
The formula (\ref{boyer_hc}) together with the isomorphism (\ref{another_iso}) imply the existence of
a basis $(o, \rho)$ of $\Gamma(\spp_+)$ such that
\[
\nabla o=A\otimes o, \quad \nabla \rho=A\otimes \rho.
\]
Let  $\iota\in \Gamma(\spp_+)$
be the covariantly constant spinor defining $N$. Then $\iota= h_1 o+h_2 \rho$ for some functions
$h_1, h_2$ on $\X_\C$. But then
\[
\nabla (h_1 o+h_2\rho)=0
\]
gives $h_1=\mbox{const}\cdot h_2$ and $A=-d\ln{(h_1)}$ so $g$ is conformal to a Ricci--flat metric.
\koniec
\vskip 5pt
\begin{prop}
\label{two_divisors}
Let $g$ be an $SL(2)$--invariant cohomogeneity--one metric with anti--self--dual Weyl tensor
on  $\X_{\C}$,  such that the quartic (\ref{quartic}) does not identically vanish. Then the
quartic vanishes on each twistor line at one point to order $4$ if and only if there exists
a null--K\"ahler structure $N$ (compatible with some metric in the conformal class of $g$) 
which is Lie--derived by the $SL(2)$ action.
\end{prop}
\noindent
{\bf Proof.}
First assume that the quartic $s$ vanishes to order $4$. Therefore its pull-back (\ref{quartic1})
to ${\mathcal F}$ factories as
\be  
\label{quartic2}
q^*(s)=\Big(\sum_i \iota_i\pi^i\Big)^4,
\ee
where $[\pi]$ are homogeneous coordinates on the fibres of $\PP(\spp_+)$ and 
$\iota$ is a section of $\spp_+$. This is a pull--back from ${\mathcal Y}$, so it is constant
along the twistor distribution (\ref{tetradlax}), i.  e.
\be
\label{laxs}
l_i(q^*(s))=0, \quad i=1, 2
\ee
which implies that $\iota$ satisfies the rank--one conformally invariant twistor equation
\[
\nabla_{i(j}\iota_{k)}=0,
\]
where $\nabla$ is the spin connection  on $\spp_+$ induced by the Levi--Civita connection of
$g$. In \cite{D02} it was shown that in this case there exists a conformal rescaling 
\[
\hat{g}=F^2 g, \quad \hat{\iota}= F\iota, \quad \hat{\varepsilon}_+= 
F\varepsilon_+
\]
such that
$\hat{\iota}$ is parallel with respect to the Levi--Civita connection of $\hat{g}$.
This section defines the null--K\"ahler structure via the isomorphism (\ref{another_iso}).

Conversely, let us assume that $N$ is an $SL(2)$--invariant null--K\"ahler structure. Then,
by Theorem \ref{mytheo}, 
$N$ gives rise to a divisor line bundle over ${\mathcal Y}$ given by a canonical  section
of $\kappa^{-1/4}$. The zero--set of this section pulls back to the hypersurface
\be
\label{omega_zero}
\Omega(l_1, l_2)^{1/2}=0
\ee
in ${\mathcal F}$, 
where $\Omega$ is the fundamental two--form of $N$ given by (\ref{hermitean_form}),
and $l_1, l_2$ span the twistor distribution. As $N$ is invariant under the $SL(2)$ action, 
the holomorphic vector fields $\widetilde{R}_\alpha, \alpha=1, 2, 3$ on ${\mathcal Y}$
preserve the canonical section of $\kappa^{-1/4}$, and so are tangent to the surface
of vanishing (\ref{omega_zero}), or equivalently
\[
\sum_i\pi_i \iota^i=0,
\]
where $\iota$ is the section of $\spp_+$ corresponding, via (\ref{another_iso}), to $N$.
Every twistor line intersects this surface at one point given,  in homogeneous coordinates, by
$[\pi]=[\iota]$. Therefore the quartic $s$ (which by assumption does not vanish identically)
vanishes at this point to order 4.
\koniec

\vskip5pt
We are now going to use the structure of the twistor
distribution (\ref{tetradlax}) to establish Theorem \ref{theorem2} from the introduction.

\noindent
{\bf Proof of Theorem \ref{theorem2}}
Let $t:\X_\C\rightarrow \C$ be a function parametrising the orbits
of $G=SL(2)$ in $\X_\C$ such that 
\[
{\mathcal L}_{R_\alpha} t=0, \quad \alpha=1, 2, 3.
\]
We can choose coordinates on $G$ such that\footnote{This form is general, but 
is different from the one
usually used (see \cite{tod1, tod2, tod3, Hitchinis}), where the vector field $\p/\p t$ is not null, and orthogonal to the orbits of $G$. We shall
explain the connection between the two forms (which are equivalent) 
in \S\ref{sectionequiv}.}
\[
g=\sum_{\alpha, \beta}\gamma_{\alpha\beta}(t) \sigma^{\alpha}\otimes\sigma^{\beta}
+\sum_{\alpha} n_\alpha(t) (\sigma^{\alpha} \otimes dt+dt \otimes \sigma^{\alpha}),
\]
where $\gamma$ is a symmetric 3 by 3 matrix and $n$ is a vector with components depending on $t$.

We can express the frame $E_{ij}$ in the distribution (\ref{tetradlax}) in 
terms of the vector field $\p_t$, and  three linearly independent vector fields $P, Q, R$ tangent to $G$ which are $t$--dependent and invariant under left translations. A convenient choice
which gives rise to the general metric of the form (\ref{metric_form})
is
\be \label{tetrad}  
E_{11} = Q,\quad 
E_{22} = P, \quad 
E_{12}= -2\frac{\p}{\p t},
\quad E_{21}=2\frac{\p}{\p t}-R.
\ee
The invariance condition implies that the functions 
$f_1$ and $f_2$ in (\ref{tetradlax}) are constant on $G,$ and so depend only on
$\lambda$ and  $t$. The quartic $s$ is proportional to 
$(f_1-\lambda f_2)$. By Proposition (\ref{two_divisors}) it must have a 
quadrupole zero which can be 
moved to $\lambda=\infty$ by a M\"obius transformation. Using the freedom
in the frame rotations of the frame we set
$(f_1=-1, f_2=0)$ so that
\[
l_1=Q+2\lambda\frac{\p}{\p t}-\frac{\p}{\p \lambda}, \quad
l_2=2\frac{\p}{\p t}-R-\lambda P.
\] 
Now consider a pair of linear combinations of $l_1$ and $l_2$ (\ref{tetradlax}) given 
by
\begin{eqnarray}
\label{tetradisolaxL}
m_1 &:=& \frac{l_1 -\lambda l_2}{f_1 + \lambda f_2} =  \frac{\p}{\p \lambda} - Q-
\lambda R-\lambda^2 P,\\ \nonumber 
m_2 &:=& \frac{1}{2}\frac{f_2 l_1 - f_1 l_2}{f_1 + \lambda f_2} = \frac{\p}{\p t} -\frac{1}{2}R-\frac{1}{2}\lambda P.
\end{eqnarray}
Since the conformal class is anti--self--dual, the integrability condition 
(\ref{intlm}) 
implies that $[m_1, m_2] =0,$ modulo $m_1$ and $m_2.$ As the Lie bracket $[m_1, m_2]$ does not
contain $\p_\lambda$ or $\p_t,$ it must be identically zero which 
yields
\be
\label{isom_airy}
\dot{P}=0, \quad \dot{Q}=\frac{1}{2}[R, Q]+\frac{1}{2}P, \quad \dot{R}=\frac{1}{2}[P, Q],
\ee
where $\dot{P}=dP/dt$ etc.

The system (\ref{isom_airy}) underlies  the isomonodromic problem with irregular singularity
of order four. To make this transparent we shall use the representation of $\mathfrak{sl}(2)$ by
2 by 2 matrices rather than vector fields, and make the replacements
\[
L_1\rightarrow \left(\begin{array}{cc}
1/2&0\\
0&-1/2
\end{array}
\right),\quad 
L_2\rightarrow \left(\begin{array}{cc}
0&1\\
0&0
\end{array}
\right),\quad
L_3\rightarrow \left(\begin{array}{cc}
0&0\\
1&0
\end{array}
\right).
\]
The associated the Lax pair (\ref{tetradisolaxL}) is the isomonodromic Lax pair for Painlev\'e II if $P$
is diagonalisable, and is gauge equivalent to the isomonodromic Lax pair for Painlev\'e I or leads to a solvable equation if
$P$ is nilpotent. The system
(\ref{isom_airy}) underlies  the isomonodromic problem with irregular singularity
of order four. To set this problem up
consider a $2\times 2$ matrix
\[
{\mathcal A}(t, \lambda)=Q+ \lambda R+\lambda^2 P, 
\]
where $\lambda\in\CP^1$, and $P, Q, R$ are elements of 
a matrix Lie algebra $\mathfrak{g}=\mathfrak{sl}(2)$  which also depend
on a parameter $t$.  When $t$ is allowed to vary on the complex plane,
the matrix fundamental solution $\Psi$  of the ODE 
\[
\frac{d \Psi}{d \lambda}={\mathcal A} \Psi
\]
depends on $\lambda$ and $t$.
The monodromy  around the pole of order four  at $\lambda=\infty$
does not depend on $t$ if $\Psi$ satisfies
\cite{Jimbo}
\be
\label{over_det}
\frac{\p \Psi}{\p \lambda}-{\mathcal A}\Psi=0, \quad
\frac{\p \Psi}{\p t}-{\mathcal B}\Psi=0,\quad
\mbox{where}\quad {\mathcal B}:=\frac{1}{2}R+\frac{1}{2}\lambda P.
\ee
The compatibility conditions
\be
\label{comcon}
\p_t {\mathcal A}-\p_{\lambda}{\mathcal B}+[{\mathcal A}, {\mathcal B}]=0
\ee
for the overdetermined linear system (\ref{over_det}) 
reduce to system of  nonlinear matrix ODEs
(\ref{isom_airy})  for $(P, Q, R)$. We shall follow the seminal work of 
\cite{mason0, mason} - but make different gauge choices - to reduce this system further to a single ODE. There are
three gauge equivalence classes to consider. The first two lead to Painlev\'e ODEs
and the last one is completely solvable.
\begin{itemize}
\item
If $P$ is diagonalisable, and 
$\mathfrak{g}=\mathfrak{sl}(2, \C)$, then  without loss of generality we can 
take \cite{Jimbo, mason} 
\[
P=2L_1, \quad R= uL_2-2\frac{z}{u} L_3, \quad
Q=(2z+t)L_1-uy L_2 -\frac{2yz+\frac{1}{2}-\alpha}{u} L_3,
\]
where $u, y$ and $z$ are functions of $t$.  Equations 
(\ref{isom_airy}) become
\[
\dot{u}=-yu, \quad \dot{z}=-2yz+\alpha-\frac{1}{2}, \quad\dot{y}=z+y^2+
\frac{t}{2},
\]
which imply
\be
\ddot{y}=2y^3+ty+\alpha,
\ee
where $\alpha$ is a constant parameter. This is 
the Painlev\'e II equation.
\item  If $P$ is nilpotent, then (as it is also constant) we can chose it to be
$L_2$. Assume that $\mbox{Tr}(PR)\neq 0$, and perform a gauge transformation
\[
{\mathcal A}\rightarrow \gamma {\mathcal A}\gamma^{-1}+\p_\lambda\gamma
\cdot \gamma^{-1}, \quad {\mathcal B}\rightarrow \gamma {\mathcal B}
\gamma^{-1}+\p_t\gamma
\cdot \gamma^{-1}
\]
with the group element $\gamma=\gamma(t)$ such that
\[
\p_t \gamma\cdot \gamma^{-1}= yP, \quad\mbox{where}\quad
y\equiv\frac{1}{8}\mbox{Tr}(R^2).
\]
Then
\[
P=L_2, \quad R= yL_2+4L_3, \quad Q=-2z L_1+\Big(y^2+\frac{t}{2}\Big)L_2-4yL_3
\]
and
\[
{\mathcal B}=\frac{\p}{\p t}-\frac{1}{2}(R+yL_2)-\frac{1}{2}\lambda P.
\]
The compatibility conditions (\ref{comcon}) give
\be
\label{PI}
\dot{y}=z, \quad \dot{z}=6y^2+t\quad\mbox{so that}\quad \ddot{y}=6y^2+t
\ee
which is the Painlev\'e I equation.
\item Finally let us consider the case where $P$ is nilpotent, and 
$\mbox{Tr}(PR)$=0. We shall
set $P=L_2$, as in the case leading to (\ref{PI}), and consider 
\[
R=\sum_{\alpha} r_{\alpha}(t) L_{\alpha}, \quad Q=\sum_{\alpha} q_{\alpha}(t) L_{\alpha}, \quad r_3=\mbox{Tr}(PR)=0
\]
for some functions $r_\alpha, q_\alpha, \alpha=1, 2, 3$ of $t$. 
The 3rd equation in (\ref{isom_airy}) gives
\[
q_1=-2\dot{r}_2, \quad  \dot{r}_1=q_3.
\]
The 2nd equation in (\ref{isom_airy}) gives
\be
\label{casel3}
2\ddot{r}_1+\dot{r}_1r_1=0, \quad 2\ddot{r}_2+\dot{r}_1r_2=0, \quad
2\dot{q}_2-2r_2\dot{r}_2-q_2r_1-1=0.
\ee
The first of these equations has a singular solution $r_1=\mbox{const}$ which eventually leads to a degenerate tetrad (\ref{tetrad}). We therefore focus on the regular solution, and absorb
two constants of integrations  in $r_1$ into an affine transformation of  $t$ which results in a constant
rescaling of the metric. The remaining two equations can also be solved:
\begin{eqnarray*}
r_1&=&4\tanh{(t)},  \quad
r_2=(a+bt)r_1-4b,\\
q_2&=&\frac{1}{4}\sinh{(2t)}-\frac{d}{dt}\Big((a+bt)r_2\Big)+c\cdot\cosh{(t)}^2,
\end{eqnarray*}
where $a, b, c$ are the remaining constants of integration. 
\end{itemize}
Now we shall construct the conformal classes corresponding to 
Painlev\'e I and Painlev\'e II equations, and in each case find a null--K\"ahler metric in the conformal class. These structures will be expressed in terms of left--invariant 
one--forms (\ref{one_formss}).
Each conformal class is represented by a covariant metric dual to 
(\ref{tetrad_0})
\[
\frac{1}{2}(e^{11}\otimes e^{22}+e^{22}\otimes e^{11}-  e^{12}\otimes e^{21}-e^{21}\otimes e^{12}).
\]
The null--K\"ahler two--form $\Omega$ can be read off from the divisor quartic
(\ref{quartic}). In the spinor--form
\[
\Omega=\iota\otimes\iota\otimes \varepsilon_-,
\]
where $\iota\in\Gamma({\spp_+}^*)$ is the parallel spinor. When regarded
as a section of $\PP({\spp_+}^*)$ it gives a point on $\CP^1$ which is the quadruple root of the quartic (\ref{quartic}). In our case this gives $\Omega$ proportional to $e^ {11}\wedge e^{21}$. The proportionality factor will be found together
with the conformal factor for the metric which makes $\Omega$ parallel.
\subsection{} For Painlev\'e II the one--forms dual to the tetrad (\ref{tetrad})
are
\begin{eqnarray*}
e^{22}&=&\frac{1}{2} \sigma^1+\frac{z(2z+t)}{u(4yz+1-2\alpha)}\sigma^2
+\frac{u\Big(z+\frac{t}{2}\Big)}{4yz+1-2\alpha} \sigma^3,\\
e^{11}&=&-\frac{2z}{u(4yz+1-2\alpha)}\sigma^2- 
\frac{u}{4yz+1-2\alpha}\sigma^3,\\
e^{21}&=&-\frac{2yz+1-2\alpha}{u(4yz+1-2\alpha)}\sigma^2+
\frac{yu}{4yz+1-2\alpha}\sigma^3,\\
e^{12}&=&-\frac{1}{2}dt-\frac{2yz+1-2\alpha}{u(4yz+1-2\alpha)}\sigma^2+
\frac{yu}{4yz+1-2\alpha}\sigma^3.
\end{eqnarray*}
The unique conformal factor which makes the null--K\"ahler two--from parallel
is  
\[
k=4yz+1-2\alpha
\]
so that
\be
\label{piiform1}
\Omega=2\sigma^3\wedge\sigma^2
\ee
and the metric is given by (\ref{metric_form}) with
\begin{eqnarray}
\label{piiform}
\gamma&=&\left(\begin{array}{ccc}
0&\frac{z}{u}&\frac{u}{2}\\
\frac{z}{u}&  \frac{8z^3+(8y^2+4t)z^2+(8-16\alpha)yz+8(\alpha-\frac{1}{2})^3}{ku^2} 
&-2\frac{2y^2z+(1-2\alpha)y-z(2z+t)}{k}  \\
\frac{u}{2}&-2\frac{2y^2z+(1-2\alpha)y-z(2z+t)}{k}&\frac{u^2(2y^2+2z+t)}{k}
\end{array}
\right),\nonumber\\
n&=&(0, \frac{2yz-2\alpha+1}{2u}, -\frac{yu}{2}).
\end{eqnarray}
The Weyl tensor is ASD if Painlev\'e II holds. 
\subsection{}
In the case of Painlev\'e I, we first read--off the tetrad
form the form of ${\mathcal A}$ and ${\mathcal B}$ to be
\[
E_{11}=Q, \quad E_{22}=P, \quad E_{12}=-2\frac{\p}{\p t}+yL_2, \quad
E_{21}=2\frac{\p}{\p t}-(R+yL_2).
\]
Computing the dual tetrad gives
\begin{eqnarray*}
e^{22}&=&\frac{y^2+\frac{t}{4}}{z}\sigma^1+\sigma^2-\frac{y}{4}\sigma^3+\frac{y}{2}dt,\\
e^{11}&=&-\frac{1}{2z}\sigma^1,\\
e^{21}&=&\frac{y}{2z}\sigma^1-\frac{1}{4}\sigma^3,\\
e^{12}&=&\frac{y}{2z}\sigma^1-\frac{1}{4}\sigma^3-\frac{1}{2}dt
\end{eqnarray*}
and rescalling the resulting metric and two--form by $k=16z(t)$ gives
\be
\label{piform1}
\Omega=2\sigma^3\wedge\sigma^1.
\ee
This is the only scaling factor which makes $\Omega$ closed. Using the same conformal factor for the metric yields
and $g$ in the form (\ref{metric_form}) with
\be
\label{piform}
\gamma=-\left(\begin{array}{ccc}
\frac{12y^2+2t}{z}&4&-3y\\
4&0&0\\
-3y&0&z
\end{array}
\right),\quad n=(0, 0, -z).
\ee
The two--form $\Omega$ is now parallel, as
\[
\nabla\Omega=\Big(\frac{6y^2+t-\dot{z}}{z}\sigma^1+
\frac{3(z-\dot{y})}{z}\sigma^2\Big)\otimes\Omega=0
\]
where we used (\ref{PI}). The null--K\"ahler structure is given 
by
\[
N=\frac{1}{2}\Big(\sigma^3\otimes L_2-\frac{2}{z}\sigma^1\otimes\frac{\p}{\p t}\Big).
\]
Computing the self--dual part of the Weyl tensor we find that it vanishes as a consequence of Painlev\'e I
\subsection{}
Computing the dual tetrad in the  case (\ref{casel3})
gives
\begin{eqnarray*}
e^{22}&=& (b\coth{(t)}-(a+bt))\sigma^1+\sigma^2\\
&&+
\Big(2b(a+bt)\coth{(t)}+b^2 \cosh{(t)}^2
-\frac{1}{8}\sinh{(t)}\cosh{(t)}^3-(a+bt)^2-\frac{c}{4}\cosh{(t)}^2  \Big)\sigma^3 ,\\
e^{11}&=&\frac{1}{4}\cosh{(t)}^2\sigma^3 ,\\
e^{21}&=&-\frac{1}{4}
\coth{(t)}\sigma^1-\frac{1}{2}\Big(b\cosh{(t)}^2+(a+bt)\coth{(t)}\Big)\sigma^3 ,\\
e^{12}&=&-\frac{1}{2}dt -\frac{1}{4}\coth{(t)}\sigma^1-\frac{1}{2}\Big(b\cosh{(t)}^2+(a+bt)\coth{(t)}\Big)\sigma^3
\end{eqnarray*}
and rescaling the resulting metric and two--form by $k=8\sinh{(t)}/\cosh{(t)}^3$ gives
\[
\Omega=\sigma^3\wedge\sigma^1.
\]
This is the only scaling factor which makes $\Omega$ closed. Using the same conformal factor for the metric yields
$g$ in the form (\ref{metric_form}) with
\begin{eqnarray}
\label{solv_form}
\gamma&=&\left(\begin{array}{ccc}
\frac{2}{\sinh{2t}} &0&(a+bt)\coth{(2t)}\\
0&0&-\tanh{(2t)}\\
(a+bt)\coth{(2t)} & -\tanh{(2t)} & 4(a+bt)^2\coth{(t)}+\frac{1}{8}\sinh{(2t)}^2(1+2c\cdot 
\coth{(t)})
\end{array}
\right),\nonumber\\
n&=&(\cosh{(t)}^{-2}, 0, 2(a+bt)\cosh{(t)}^{-2}+2b\tanh{(t)}).
\end{eqnarray}
The two--form $\Omega$ is parallel and  the Weyl tensor is anti--self--dual. The Ricci--tensor is 
\[
r=\frac{1}{2}\sinh{(t)}\cosh{(t)}^3\sigma^3\otimes\sigma^3.
\]
\koniec
{\bf Remarks}
\subsection{} All solutions to PI are transcendental, and PII admits 
solutions expressible in terms of know functions only for integer and 
half--integer values of the parameter $\alpha$. These will lead to explicit metrics. On the other hand all metrics 
arising from (\ref{solv_form}) are explicit. The simplest 
is obtained by setting $a=b=c=0$. Setting
$\tau=\tanh{(t)}$ gives
\[
g=\frac{1-\tau^2}{\tau}\sigma^1\otimes\sigma^1 +\sigma^1\otimes d\tau+d\tau\otimes\sigma^1-
2\tau(\sigma^2\otimes\sigma^3+\sigma^3\otimes\sigma^2)+
\frac{\tau^2}{2(1-\tau^2)^2}\sigma^3\otimes\sigma^3.
\]
\subsection{} The null--Kahler structures
arising from PI and and PII can be distinguished  by examining the restriction of kernel
of the endomorphism $N$ to the orbits of $SL(2)$. In the PI case, this kernel - when regarded as the element of the Lie algebra of $SL(2)$ is 
nilpotent, but in the PII case it is not. This can be seen directly from
(\ref{piform1}) and (\ref{piiform1}).
The kernel of $N$ spans a 2-parameter of $\alpha$--surfaces in $\X_\C$ which, at the twistor level, corresponds to a divisor hypersurface 
${\mathcal N}\subset {\mathcal Y}$ which meets each twistor line to 
order 4. This hypersurface is preserved by the 
$SL(2, \C)$ action on the twistor space,  
so $\exists \tau\in\mathfrak{sl}(2)$ s. t.
$\phi(\tau)=0$ where $\phi$ is the
vector bundle homomorphism \cite{Hitchinis}
\[
\phi:\mathfrak{sl}(2, \C)\times{\mathcal Y}\rightarrow T{\mathcal Y}.
\]
The element $\tau$ is nilpotent in the case of Painlev\'e I and the solvable example, and semisimple in the case of Painlev\'e II.
\subsection{} 
\label{sectionequiv}
The usual form of cohomogeneity--one metrics \cite{tod1, tod2, tod3, Hitchinis}
is
\[
g=\frac{1}{4}dt^2+\sum_{\alpha, \beta}h_{\alpha\beta}(t)(\sigma^\alpha\otimes\sigma^\beta+\sigma^\beta\otimes \sigma^\alpha).
\]
This arises from a frame of the form
\[
E_{11}=Q,\quad E_{22}=P,\quad E_{12}=-2\frac{\p}{\p t}-\frac{1}{2}R,\quad E_{21}=2\frac{\p}{\p t}-\frac{1}{2}R.
\]
Following the argument above which lead to the Lax pair (\ref{tetradisolaxL})
gives the system of ODEs
\cite{CD}
\[
\dot{Q}=\frac{1}{4}[R, Q]+\frac{1}{2}P, \quad \dot{R}=\frac{1}{2}[P, Q], \quad
\dot{P}=\frac{1}{4}[P, R]
\]
together with the isomonodromic Lax pair $[m_1, m_2]=0$, where
\be
\label{system_new}
m_1=-l_1+\lambda l_2=\p_\lambda-{\mathcal A}, \quad m_2=\frac{1}{2}l_2=\p_t-{\mathcal B}
\ee
where now
\[
{\mathcal A}=Q+\lambda R+\lambda^2 P,  \quad {\mathcal B}=\frac{1}{4}R+\frac{1}{2}\lambda P.
\]
The gauge transformation
\[
{\mathcal A}\rightarrow \gamma {\mathcal A} \gamma^{-1}+\p_\lambda\gamma\cdot\gamma^{-1},
\quad {\mathcal B}\rightarrow \gamma {\mathcal B} \gamma^{-1}+\p_t\gamma\cdot\gamma^{-1}
\]
with $\gamma=\gamma(t)$ such that $\gamma^{-1}\cdot\dot{\gamma}=\frac{1}{4}R$ brings this Lax pair to the form (\ref{tetradisolaxL}), and the system (\ref{system_new}) to the form (\ref{isom_airy}).  Indeed, we can verify that $\widetilde{P}\equiv \gamma \cdot P\cdot\gamma^{-1}$ is constant,
and the other two equations also hold  with   $\widetilde{Q}\equiv \gamma \cdot 
Q\cdot\gamma^{-1}$ and $\widetilde{R}\equiv \gamma \cdot R\cdot\gamma^{-1}$.
Thus the two forms of the metric are equivalent by a diffeomorphism.
\subsection{}
The isomonodromic Lax pair (\ref{tetradisolaxL}) arising as the combination of $(l_1, l_2)$ can be constructed invariantly using the notation  from the proof
of Proposition \ref{quartic_lemma} as follows:
Let $T$ be the vector field dual to $dt$ with respect to the isomorphism between
tangent and cotangent bundle given by the metric. This vector is normal to the orbits
of $SL(2)$, and is given by $T=\sum_{i, j} T^{ij} E_{ij}$ in the frame $E_{ij}$. First note that
\[
l_i(t)=\sum_j \pi^j T_{ij},  \quad\mbox{so that}\quad
\sum_{i, j} T^{ij}\pi_j l_i(t)=\sum_{i, j, k, m} \pi_i\pi_j\varepsilon_{km}T^{ki}T^{mj}=0
\]
by symmetry. This implies that
\[
\sum_{i, j} T^{ij}\pi_jl_i=\sum_{i, j} T^{ij}\pi_j f_i\frac{\p}{\p\lambda}+\sum_{i, j, k}
T^{ij}\pi_j\pi^k E_{ik},
\]
where the second term on the RHS is tangent to the orbits of $SL(2)$, so is in the span of the left--invariant vector fields and does not contain $\p/\p t$. Using (\ref{new_form}) we identify the  multiple
of $\p/\p \lambda$ as the quartic $s$. Assuming that this quartic is not identically zero we define
an $\OO(-2)$--valued vector field
\[
m_1=\frac{\sum_{i, j}T^{ij}\pi_j l_i}{s}.
\]
The second vector field $m_2$ does not contain the $\p/\p\lambda$ terms, and is defined by
\[
m_2=\frac{1}{2}\frac{\sum_{i, j} \varepsilon^{ij}f_i l_j}{s}.
\]
This agrees with (\ref{tetradisolaxL}).
If the metric is given by (\ref{metric_form}), then
\[
T=-\frac{1}{\sum_{\alpha, \beta}\gamma^{\alpha\beta} n_\alpha n_\beta}\Big(\frac{\p}{\p t}-\sum_{\alpha, \beta}\gamma^{\alpha\beta}
n_\alpha L_\beta\Big),
\]
where $\gamma^{\alpha\beta}(t)$ is the inverse matrix of $\gamma_{\alpha\beta}(t)$.
\subsection{}
Selecting a one--parameter family of transformations $\R^{*}$ in $SL(2, \R)$  generated by a non--null Killing vector $K$ reduces the
null--K\"ahler ASD conditions to an Einstein--Weyl structure in 2+1
dimension which additionally admits a parallel weighted null vector field.
Such structures correspond to solutions of the dispersionless Kadomtsev-Petviashvili equation, and in \cite{DTT},  such solutions were
constructed and shown to be constant on central quadrics and expressed in terms of solutions to Painlev\'e I or Painlev\'e II.
\subsection{}
 If the Lie algebra  $\mathfrak{g}$ underlying (\ref{isom_airy}) is instead 
the Bianchi II algebra 
then the isomonodromic condition is the (derivative of ) the Airy equation
(see \cite{CD}, where a class of ASD null--K\"ahler four manifolds has been constructed).
In \cite{Hitchinis, mmw} it was instead assumed that the quartic $s$ has four distinct zeros, and that $G=SL(2, \C)$ which lead to the isomonodromic Lax pair  \cite{Jimbo} for Painlev\'e VI.
If $s$ has two double zeroes then the conformal class contains an Einstein 
metric \cite{tod2, Hitchinis}, and the isomonodromic Lax pair leads to 
Painlev\'e III.

\end{document}